\documentclass[12pt]{article}

\usepackage{amsthm}
\usepackage[centertags]{amsmath}
\usepackage{amssymb}
\usepackage{amsfonts}
\usepackage{graphics}

\usepackage{cite}

\DeclareMathOperator{\diam}{\mathrm{diam}}
\DeclareMathOperator{\Lip}{\mathrm{Lip}}

\def\M{\mathcal M}
\def\K{\mathcal K}

\def\bR{\mathbb R}
\def\la{\lambda}

\def\lan{\langle}
\def\ran{\rangle}
\def\D{\mathcal D}
\def\R{\mathcal R}
\def\E{\mathcal E}
\def\O{\mathcal O}

\def\le{\leqslant}
\def\geq{\geqslant}
\def\ge{\geqslant}
\def\wt{\widetilde}

\theoremstyle{plain}
\newtheorem{theorem}{Теорема}
\newtheorem{lemma}{Лемма}
\newtheorem{arrow}{Следствие}
\newtheorem{definition}{Определение}
\newtoks\thehProclaim
\newtheorem*{Proclaim}{\the\thehProclaim}

\newenvironment{proclaim}[1]{\thehProclaim{#1}\begin{Proclaim}}{\end{Proclaim}}

\newcommand{\bp}{\begin{proof}}
\newcommand{\ep}{\end{proof}}
\newcommand{\bl}{\begin{lemma}}
\newcommand{\el}{\end{lemma}}
\newcommand{\bt}{\begin{theorem}}
\newcommand{\et}{\end{theorem}}
\newcommand{\bd}{\begin{definition}}
\newcommand{\ed}{\end{definition}}
\newcommand{\ba}{\begin{arrow}}
\newcommand{\ea}{\end{arrow}}

\begin{document}

\thispagestyle{empty}

\vspace{.5in}
\begin{center}
\begin{minipage}{14cm}
\begin{center}
\Large \bfseries
Point distributions in compact metric spaces 
\end{center}
\end{minipage}
\\[0.5in]
\begin{tabular}{c}
{\large\bf M. M. Skriganov} \\[6pt]
St.Petersburg Department of\\
Steklov Mathematical Institute\\
Russian Academy of Sciences \\[6pt]
E-mail: {\tt maksim88138813@mail.ru}

\end{tabular}

\end{center}

\vspace{.5in}
\vspace{.5in}
 We consider finite point subsets (distributions) in compact metric spaces.
Non-trivial bounds for sums of distances between points of distributions and
for discrepancies of distributions in metric balls are given in the case of
general rectifiable metric spaces (Theorem 1.1)

 We generalize Stolarsky's invariance principle to distance-invariant spaces
(Theorem 2.1), and for arbitrary metric spaces we prove a probabilistic
invariance principle (Theorem 3.1).

 Furthermore, we construct partitions of general rectifiable compact metric
spaces into subsets of equal measure with minimum average diameter (Theorem 4.1).

\vskip1cm


\vfil

\newpage












  
	


\theoremstyle{plain}







\section*{Contents}

\noindent 1. Introduction

\noindent 2. The invariance principle for distance-invariant spaces

\noindent 3. Equal measure partitions and the probabilistic invariance principle 

\noindent 4. Construction of equal measure partitions. Proof of Theorem 1.1

\noindent    References 

\medskip

\section*{1. Introduction}

Let $\M$ be a compact separable metric space with a metric $\rho$ 
and a finite non-negative Borel measure  $\mu$ normalized by $\mu(\M)=1$. 

For any $N$-point subset (distribution) 
$\D_N=\{x_1,\dots, x_N\}\subset \M$ we put 
\begin{equation}
\rho [\D_N]=\sum\limits^N_{i,j=1}\rho(x_i,x_j)
\tag{1.1}
\end{equation} 
and denote by $\langle \rho\rangle$ the mean value of the metric $\rho$,
\begin{equation}
\langle \rho \rangle=\iint\limits_{\M\times \M}\rho(y_1,y_2)\,d\mu 
(y_1)\,d\mu(y_2).
\tag{1.2}
\end{equation}

We write $B_r(y)=\{x:\rho(x,y)\le r\}$, $r\in T$, $y\in \M$, for the ball of 
radius $r$ centered at $y$, and of volume $\mu(B_r(y))$. Here  
$T=\{r:r=\rho(y_1,y_2),y_1,y_2\in \M\}$ is the set of radii, $T\subset [0,L]$, where 
$L=\sup \{r=\rho(y_1,y_2):y_1,y_2\in \M\}$ is the diameter of $\M$.

The local discrepancy $\D_N$ is defined by 
\begin{align*}
\Lambda[B_r(y),\D_N]&=\#\{B_r(y)\cap \D_N\}-N\mu (B_r(y))
\\
&=\sum\limits_{x\in \D_N} \Lambda (B_r(y),x),
\tag{1.3}
\end{align*}  
where
\begin{equation}
\Lambda(B_r(y),x)=\chi(B_r(y),x)-\mu(B_r(y))
\tag{1.4}
\end{equation}
and $\chi(\E,x)$ is the characteristic function of a subset  $\E\subset \M$. 

We put
\begin{equation}
\la_r[\D_N]=\int\limits_{\M}\Lambda[B_r(y),\D_N]^2\,d\mu(y)
\tag{1.5}
\end{equation}
This formula can be written as
\begin{equation}
\la_r[\D_N]=\sum\limits_{y_1,y_2\in \D_N}\la_r(y_1,y_2),
\tag{1.6}
\end{equation}
where
\begin{equation}
\la_r(y_1,y_2)=\int\limits_{\M}\Lambda(B_r(y),y_1)\Lambda(B_r(y),y_2)\,d\mu(y).
\tag{1.7}
\end{equation}

Let $\xi$ be a finite non-negative measure on the set of radii $T$. 
We put
\begin{equation}
\la[\xi,\D_N]=\int\limits_T\la_r[\D_N]d\xi(r)=\sum\limits_{y_1,y_2\in 
\D_N}\la(\xi,y_1,y_2),
\tag{1.8}
\end{equation}
where
\begin{equation}
\la(\xi,y_1,y_2)=\int\limits_T\la_r(y_1,y_2)\,d\xi (r).
\tag{1.9}
\end{equation}

The quantity $\la[\xi,\D_n]^{1/2}$ is known as the $L_2$-discrepancy of a distribution  
$\D_N$ in bolls $B_r(y)$, $r\in T$, $y\in \M$ with respect to the measures $\mu$ and $\xi$.

Introduce the following extremal characteristics 
\begin{align}
\rho_N(\M)&=\sup\limits_{\D_N}\rho [\D_N],
\tag{1.10}
\\
\la_N(\xi,\M)&=\inf\limits_{\D_N}\la [\xi, \D_N],
\tag{1.11}
\end{align} 
where the supremum and infimum are taken over all $N$-point distributions $\D_N\subset \M$.

The study of the characteristics (1.10) and (1.11) is a subject of the geometry of  
distances and the discrepancy theory, see \cite{2, 5, 10}.

In the present paper we will show that non-trivial bounds for the quantities (1.10)  
and (1.11) can be obtained under very general conditions on spaces 
$\M$, metrics $\rho$ and measures $\mu$ and $\xi$.

It is convenient to introduce the concept of $d$-rectifiable spaces,
that will allow us to compare the metric and measure on $\M$ with
Euclidean metric and Lebesgue measure on  $\bR^d$. The concept of
rectifiability is well known in the geometric measure theory, see \cite{11}. 
Here this terminology is adapted to our purposes.

Recall that a map $f:\O\subset \bR^d\to \M$ is a Lipschitz map if
\begin{equation}
\rho(f(Z_1),f(Z_2))\le c\|Z_1-Z_2\|, Z_1,Z_2\in \O,
\tag{1.12}
\end{equation}
with a constant $c$, and the smallest such constant is called the Lipschitz constant of $f$ 
and denoted by $\Lip(f)$; in (1.12) $\|\cdot \|$ is Euclidean metric in $\bR^d$. 

\begin{proclaim}{Definition 1.1} A compact metric space $\M$ with a metric $\rho$ and 
a measure $\mu$ is called $d$-{\it rectifiable} if there exist a measure $\nu$ on 
the $d$-dimensional unite cube $I^d=[0,1]^d$ absolutely continuous with respect 
to Lebesgue measure, a measurable subset $\O\subset I^d$, and an injective Lipschitz 
map $f:\O\to \M$,  such that

{\rm(i)} $\mu(\M\setminus f(\O))=0$,

{\rm(ii)} $\mu(\E)=\nu(f^{-1}(\E\cap f(\O))$ for any $\mu$-measurable 
subset $\E\subset \M$.
\end{proclaim}

Since the map $f$ is injective, the formula
\begin{equation}
\nu(\K\cap \O)=\mu(f(\K\cap \O))
\tag{1.13}
\end{equation}
is well defined for any measurable subset $\K\subset I^d$. Also, we can
assume that the measure $\nu$ is concentrated on $\O$ and $\nu(\O)=\mu(f(\O))=\mu(\M)=1$.

Simple examples of $d$-rectifiable spaces can be easily given. Any
smooth (or piece-wise smooth) compact $d$-dimensional manifold is 
$d$-rectifiable, if in the local coordinates the metric satisfies (1.12),
and the measure is absolutely continuous with respect to Lebesgue measure. 
Particularly, any compact $d$-dimensional Riemannian manifold with Riemannian metric 
and measure is $d$-rectifiable. In this case, it is known that the condition (1.12) 
is true, see ~\cite[Chapter~I, Proposition 9.10]{8}, while the condition on the measure
is obvious because the metric tensor is continuous. We refer to \cite{11} for much more
exotic examples of rectifiable spaces.

In the present paper we will prove the following theorem.

\begin{proclaim}{Theorem 1.1} Let $\M$ be a compact $d$ --rectifiable  
metric space. Then for each $N$ we have 
\begin{equation}
\rho_N(\M)\geq N^2\lan \rho \ran -d2^{d-1}\Lip(f)N^{1-\frac{1}{d}}.
\tag{1.14}
\end{equation}

If additionally, the measure $\xi$ on the set of radii $T$ satisfies
the condition
$$
\xi([a,b))\le c_0(\xi)|a-b|, \quad a\le b, a, b\in T,
$$
with a constant $c_0(\xi)>0$, then for each $N$ we also have
\begin{equation}
\la_N(\xi,\M)\le d2^{d-2}\Lip(f)c_0(\xi)N^{1-\frac{1}{d}}.
\tag{1.15}
\end{equation}

In {\rm(1.14)} and {\rm(1.15)} $\Lip(f)$ is the Lipschitz constant of 
the map $f$ in the definition of $d$-rectifiability of the space $\M$.
\end{proclaim} 

Under such general assumptions one can not expect that the bounds (1.14) 
and (1.15) are the best possible, and one can give examples of  
$d$-dimensional manifolds where the bounds (1.14) and (1.15) can be 
improved. At the same time, it is known that for the $d$-dimensional spheres $S^d\subset 
\bR^{d+1}$ with the rotation invariant metric and measure  
the bounds (1.14) and (1.15) are the best possible.

For spheres $S^d$ bounds of the type (1.14) and (1.15) were established by Stolarsky 
\cite{14}. The opposite bounds 
\begin{equation}
\rho_N(\M)\le N^2\lan \rho \ran-c_4N^{1-\frac{1}{d}}
\tag{1.16}
\end{equation}
and
\begin{equation}
\la_N(\M)\ge c_2N^{1-\frac{1}{d}}
\tag{1.17}
\end{equation}
with positive constants $c_1$ and $c_2$ independent of $N$, in the case of 
$\M=S^d$, were proved by Beck \cite{4}. We refer to \cite{2, 5, 10} for a discussion
of these results.

Spheres as homogeneous spaces $S^d=SO(d+1)/SO(d)$ are the simplest examples of
compact Riemannian symmetric spaces of rank one. All such spaces are known, see,
for example,  \cite{8}: besides the spheres they are the real, complex, 
and quaternionic projective spaces and the octonionic projective plane. 

By Theorem~1.1 the bounds (1.14) and (1.15) hold for all these spaces. 
It turns out that the opposite bounds (1.16) and (1.17) are also true for
all such spaces. This theorem will be proved in our forthcoming paper \cite{13}
with the help of methods of harmonic analysis on homogeneous spaces.
Hence, relying on Theorem~1.1 in conjunction with the results of \cite{13}, 
we are able to establish the exact order of the extremal characteristics (1.10) and (1.11) 
as $N\to \infty$ for all compact Riemannian symmetric spaces of rank one.

In the present paper we use quite elementary methods. It should be recorded that these
methods are related to the papers Alexander \cite{1} and Stolarsky \cite{14} dedicated 
to point distributions on spheres.

In Section 2 we give a generalization of Stolarsky's invariance principle to
distance-invariant spaces (Theorem~2.1).

In Section 3 we give a probabilistic version of the invariance principle
suitable for arbitrary compact metric spaces (Theorem~3.1). Using
the probabilistic invariance principle, we obtain in Section 3 the basic bounds
for the characteristics (1.10) and (1.11) in terms of equal measure
partitions of a metric space (Theorem~3.2). 
  
In Section 4 for $d$-rectifiable compact metric spaces we give an explicit construction of
equal measure partitions with the optimum order of the average diameter of the subsets
of the partition (Theorem~4.1). Using such partitions, we complete the proof of 
Theorem~1.1 in Section 4.

\section*{2.The invariance principle for distance-invariant spaces }

On the space $\M$ we introduce the following
metrics associated with the original metric $\rho$ and measure $\mu$
\begin{equation}
\rho^*(\xi,y_1,y_2)=\int\limits_T\rho^*_r(y_1,y_2)\,d\xi (r),
\tag{2.1}
\end{equation}
where 
\begin{equation}
\rho^*_r(y_1,y_2)=\mu(B_r(y_1)\Delta B_r(y_2))
\tag{2.2}
\end{equation}
is the "symmetric difference" metric for the balls 
\begin{equation}
B_r(y_1)\Delta B_r(y_2)=B_r(y_1)\cup B_r(y_2)\setminus B_2(y_1)\cap 
B_r(y_2).
\tag{2.3}
\end{equation}
Therefore, 
\begin{align}
\rho^*_r(y_1,y_2)&=\int\limits_{\M}\chi(B_r(y_1)\Delta 
B_r(y_2),y)\,d\mu(y)
\notag\\
&=\int\limits_{\M}[\chi (B_r(y_1),y)+\chi (B_r(y_2),y)
-2\chi 
(B_r(y_1),y)\chi(B_r(y_2),y)]\, d\mu(y)
\notag\\
&=\int\limits_{\M}|\chi(B_r(y_1),y)-\chi (B_r(y_2),y)|d\mu (y).
\tag{2.4}
\end{align}

For the average values of the metrics $\rho^*(\xi)$ and $\rho^*_r$ we obtain 
\begin{align}
\lan \rho^*(\xi)\ran &=\int\limits_T\lan \rho^*_r\ran d\xi(r),
\tag{2.5}  
\notag\\
\lan \rho^*_r\ran &\!=\!\iint\limits_{\M\times 
\M}\rho^*_r(y_1,y_2)\,d\mu(y_1)d\mu (y_2)  
\!\!=\!\!2\!\int\limits_{\M}[\mu (B_r(y))\!\!-\!\!\mu (B_r(y))^2]\,d\mu (y).
\tag{2.6}
\end{align} 

In view of symmetry of the metric $\rho$, we have 
\begin{equation}
\chi(B_r(y),x)=\chi (B_r(x),y)=\chi (r-\rho(x,y)),
\tag{2.7}
\end{equation}
where $\chi(t)$, $t\in \bR$ is the characteristic function of the half-axis  
$[0,\infty)$.

\begin{proclaim}{Lemma 2.1} {\rm(i)} We have the equality 
\begin{equation*}
\rho^*(\xi,y_1,y_2)=\int\limits_{\M}|\sigma(\rho(y_1,y))-\sigma(\rho(y_2,y))|d\mu(y),
\tag{2.8}
\end{equation*}
where
\begin{equation}
\sigma(r)=\xi([r,L])=\int\limits^L_r\,d\xi(t), \quad r\in T,
\tag{2.9}
\end{equation}
and $L=\sup\{r:r\in T\}$ is the diameter of $\M$.

{\rm(ii)} If the measure $\xi$ satisfies the condition 
\begin{equation}
\xi([a,b))\le c_0(\xi)|a-b|, \ a\le b, a, b\in T,
\tag{2.10}
\end{equation}
with a constant $c_0(\xi)>0$, then we also have the inequality 
\begin{equation}
\rho^*(\xi,y_1,y_2)\le c_0(\xi)\rho(y_1,y_2).
\tag{2.11}
\end{equation}
\end{proclaim}

\begin{proof} [Proof] For short, we write $\rho(y_1,y)=\rho_1$, 
$\rho(y_2,y)=\rho_2$.

(i) Using the formulas (2.1), (2.4) and (2.7), we obtain
\begin{align}
& \rho^*(\xi,y_1,y_2) \notag
\\
&=\int\limits_{\M}d\mu(y)\int\limits_Td\xi(r)\big[\chi(r-\rho_1)+\chi(r-\rho_2)  
-2\chi(r-\rho_1)\chi(r-\rho_2)\big]
\notag \\
& =\int\limits_{\M}d\mu(y)\big[\sigma(\rho_1)  
+\sigma(\rho_2)-2\sigma(\max\{\rho_1,\rho_2\})\big]
\tag{2.12}
\end{align}

Since $\sigma$ is a non-increasing function, we have 
\begin{align}
2\sigma(\max \{\rho_1,\rho_2\})&=2\min\{\sigma(\rho_1),\sigma(\rho_2\}
\notag\\
&=\sigma(\rho_1)+\sigma(\rho_2)-|\sigma(\rho_1)-\sigma(\rho_2)|.
\tag{2.13}
\end{align}
Substituting  (2.13) to (2.12), we obtain (2.8).

(ii) Suppose that $\rho_1\le \rho_2$, then using (2.9), (2.10) and the triangle inequality
for the metric $\rho$, we obtain 
\begin{align}
&|\sigma(\rho_1)-\sigma(\rho_2)|=\xi([\rho_1,L])-\xi([\rho_2,L])=\xi([\rho_1,\rho_2))
\notag\\
&\le c_0(\rho_2-\rho_1)=c_0(\rho(y_2,y_1)-\rho(y_1,y))\le 
c_0(\xi)\rho(y_1,y_2).
\tag{2.14}
\end{align}

The similar inequality holds if $\rho_1>\rho_2$. Substituting (2.14) 
to (2.8), we obtain (2.11).

The proof of Lemma 2.1 is complete.
\end{proof}

Consider the kernel (1.7). Substituting (1.4) to (1.7), we obtain
\begin{align}
\la_r(y_1,y_2)&=\int\limits_{\M}\,d\mu(y)\big[\chi(B_r(y),y_1)\chi(B_r(y),y_2)
\notag\\
&\!-\!\mu(B_r(y))\chi
(B_r(y),y_1)
\!-\!\mu(B_2(y))\chi(B_r(y),y_2)\!+\!\mu(B_r(y))^2\big].
\tag{2.15}
\end{align}

Comparing the formulas (2.4) and (2.15), we see that 
\begin{equation}
2\la_r(y_1,y_2)+\rho^*_r(y_1,y_2)=A^{(0)}_r+A^{(1)}_r(y_1)+A^{(1)}_r(y_2),
\tag{2.16}
\end{equation}
where
\begin{equation}
A^{(0)}_r=2\int\limits_{\M}(B_r(y))^2 \, d\mu(y),
\tag{2.17}
\end{equation}
\begin{align}
A^{(1)}_r(x)&=\int\limits_{\M}d\mu(y)\big[\chi(B_r(x),y)-2\mu(B_r(y))\chi(B_r(y),x)\big]
\notag\\
&=\mu(B_r(x))-2\int\limits_{\M}\mu(B_r(y))\chi(B_r(y),x)\, d\mu(y)
\notag\\
&=\mu (B_r(x))-2\int\limits_{\M}\mu(B_r(y))\chi(B_r(x),y)\, d\mu(y),
\tag{2.18}
\end{align}
here we used the formula (2.7).

Let us consider these formulas in the following special case.
A metric space $\M$ is called {\it distance-invariant}, if for each  $r\in T$ 
the volume of ball $\mu(B_r(y))$ is independent of $y\in \M$, see 
\cite{9}.

In this case, the integrals in (2.17), (2.18) can be easily evaluated and we arrive 
at the following result.

\begin{proclaim}{Theorem 2.1} Let $\M$ be a compact 
distance-invariant metric space. Then we have the relations
\begin{align}
2\la_r(y_1,y_2)&+\rho^*_r(y_1,y_2)=\lan \rho^*_r\ran
\tag{2.19}
\notag\\
2\la(\xi,y_1,y_2)&+\rho^*(\xi,y_1,y_2)=\lan \rho^*(\xi)\ran.
\tag{2.20}
\end{align}
Particularly, for any $N$-point distribution $\D_N\subset \M$ we have the invariance principle
\begin{equation}
2\la[\xi,\D_N]+\rho^*[\xi,\D_N]=N^2\lan \rho^*(\xi)\ran.
\tag{2.21}
\end{equation} 
\end{proclaim}

\begin{proof}[Proof] For short, we write $v_r=\mu(B_r(y))$. By definition, 
$v_r$ is a constant independent of $y\in \M$. and the formulas (2.17), (2.18) 
take the form
$$
A^{(0)}_r=2v^2_r,\quad A^{(1)}_r(x)=v_r-2v^2_r, \quad x\in \M.
$$

Therefore, the right side of the equality (2.16) is equal to $2(v_r-v^2_r)$.  
From the other hand, the average value (2.6) is also equal to $2(v_r-v^2_r)$. 
This proves the equality (2.19).

Integrating the equality (2.19) over $r\in T$ with the measure $\xi$, we obtain (2.20). 
Summing the equality (2.20) over $y_1$, $y_2\in \D_N$, we obtain (2.21).
\end{proof}

The typical examples of distance-invariant spaces are (finite or infinite)
homogeneous spaces $\M=G/H$, where $G$ is a compact group, $H\vartriangleleft G$ is
a closed subgroup, while $\rho$ and $\mu$  are some $G$-invariant  metric and measure on $\M$.

Numerous examples of distance-invariant spaces are known in algebraic combinatorics 
as distance-regular graphs and metric association schemes (on finite or infinite sets). 
Such spaces are characterizing even a stronger condition: the volume of the intersection of
any two balls $\mu(B_r,(y_1)\cap B_{r_2}(y_2))$ depends only on  
$r_1$, $r_2$ и $r_3=\rho(y_1,y_2)$, see \cite{3, 9}.

For spheres $S^d$ the identity (2.21) was established by Stolarsky \cite{14} and 
called the invariance principle. The original proof in \cite{14} was rather difficult, it was 
simplified in the recent papers Bilyk \cite{6} and Brauchard and Dick \cite{7}. 

Theorem 1.1 is a generalization of the invariance principle to arbitrary compact
distance-invariant spaces. Probably, the above arguments provide the most adequate 
explanation of such relations.

Notice that the formula (2.8) enables us to calculate the metric 
$\rho^*$ explicitly for some special spaces $\M$. For spheres $S^d\subset 
\bR^{d+1}$ and a special measur $\xi$ the metric $\rho^*$ is equivalent to
Euclidean metric in $\bR^{d+1}$, this fact was established in \cite{14}, see also \cite{6, 7}. 
In \cite{13} we will show that for projective spaces and specific measure $\xi$ 
the metric $\rho^*$ is equivalent to the Fubini-Study metric.

\section*{3.Equal measure partitions and the probabilistic invariance principle} 

Whether it is possible to generalize the relations (2.19), (2.20) и (2.21) to arbitrary
compact metric spaces ? At first glance the answer should be negative. Nevertheless, 
a probabilistic generalization of these relations turns out to be possible. 

First of all, we introduce some definitions and notations. Consider a partition 
$\R_{N}=\{V_i\}^N_1$ of a compact space  $\M$ into $N$ measurable subsets  $V_i \subset \M$, 
\begin{equation}
\mu (\M\setminus \bigcup\limits^{N}_{i=1} V_i)=0, \quad
\mu (V_i\cap V_j)=0,  \quad i\ne j
\tag{3.1}
\end{equation}

We write $\diam (\rho, V)=\sup \{ \rho (y_1,y_2), y_1,y_2\in V\}$
for the diameter of a subset $V\subset \M$ with respect to the metric $\rho$.
For a partition $\R_N$ we introduce  \emph{the average diameter} $\| \R_N\|_1$ by
\begin{equation}
\| \R_N\|_1=\frac{1}{N}  \sum\limits^{N}_{i=1} \diam \, (\rho, V_i)
\tag{3.2}
\end{equation}
and \emph{the maximum diameter} $\| \R_N\|_{\infty}$ by
\begin{equation}
\| \R_N\|_{\infty}=\max\limits_{1\le i\le N} \diam \, (\rho, V_i). 
\tag{3.3}
\end{equation}

A partition $\R_N= \{ V_i\}^N_1$ is \emph{an equal measure partition} if all subsets  
 $V_i$ have equal measure, $\mu (V_i)=N^{-1}$, $1\le i\le N$. 

Let an equal measure partition $\R_N=\{V_i\}^N_1$ of the space $\M$ be given. 
We introduce a probability space  $\Omega_N$ by 
\begin{equation}
\Omega_N=\prod\limits^{N}_{i=1} V_i= \{ X_N=
(x_1,\dots,x_N):x_i\in V_i , 1\le i\le N\} 
\tag{3.4}
\end{equation}
with a probability measure $\omega_N=\prod\limits^{N}_{i=1} \wt\mu_i$, where
$\wt\mu_i=N\mu|V_i$, and $\mu|V_i$ denotes the restriction of the measure $\mu$ to 
a subset  $V_i\subset \M$.

We write $\mathbb E_NF[\,\cdot\,]$ for the expectation of a random variable
$F[X_N]$, $X_N\in \Omega_N$, 
\begin{equation}
\mathbb E_NF[\, . \,] = \int\limits_{\Omega_N} 
F[X_N] \ d\omega_N= N^N \iint\limits_{V_1\times \dots \times V_N} 
F(x_1,\dots,x_N) \, d\mu (x_1)\dots \, d\mu (x_N).
\tag{3.5}
\end{equation} 

\begin{proclaim}{Lemma 3.1} Let $F^{(1)}[X_N]$ and $F^{(2)}[X_N]$,
$X_N=(x_1,\dots, x_N)\in \Omega_N$, be the following random variables 
\begin{equation}
F^{(1)}[X_N]=\sum\limits_i f (x_i), \quad
F^{(2)}[X_N] = \sum\limits_{i\ne j} f (x_i,x_j), 
\tag{3.6}
\end{equation} 
where $f(y)$ and $f(y_1,y_2)$ are integrable functions on $\M$ and $\M \times \M$, 
correspondingly. Then
\begin{align}
 \mathbb E_NF^{(1)} [\, . \,] & =N\int\limits_{\M} f (y) \, d\mu (y)
\tag{3.7}
\\
 \mathbb E_NF^{(2)} [\, . \,] &  =N^2\iint\limits_{\M\times \M} f (y_1,y_2) \, d\mu (y_1)
\, d\mu (y_2) \notag
\\
& - N^2\sum\limits^{N}_{i=1} \iint\limits_{V_i\times V_i}
f(y_1,y_2) \, d\mu (y_1)\, d\mu (y_2).
\tag{3.8}
\end{align}
\end{proclaim}                      

\begin{proof}[Proof] Substituting  (3.6) to (3.5), we obtain
$$
\mathbb E_NF^{(1)} [\, . \,] = N \sum\limits_i \int\limits_{V_i} f(y) 
\, d\mu (y) =N \int\limits_{\M} f(y)\, d\mu (y). 
$$                  
This proves (3.7). 

Substituting (3.6) to (3.5), we obtain
\begin{align*}
&\mathbb E_NF^{(2)}[\, \cdot \,] = N^2 \sum\limits_{i\ne j}
\, \iint\limits_{V_i\times V_j} f(y_1,y_2) \, d\mu (y_1)\, d\mu (y_2)
\\
&=N^2\sum\limits_{i,j} \,\iint\limits_{V_i\times V_j} f(y_1,y_2) \, d\mu 
(y_1)\, d\mu (y_2) 
-N^2\sum\limits_{i} \iint\limits_{V_i\times V_i} f(y_1,y_2) \, d\mu 
(y_1)\, d\mu (y_2) 
\\
&=N^2\iint\limits_{\M \times \M} f(y_1,y_2) \, d\mu 
(y_1)\, d\mu (y_2) - 
N^2\sum\limits_i \, \iint\limits_{V_i\times V_i}
f(y_1,y_2) d\mu (y_1)\, d\mu  (y_2). 
\end{align*}
This proves (3.8).
\end{proof}

Elements $X_N\in \Omega_N$ can be thought of as specific 
$N$-point distributions in the space $\M$ and the corresponding sums of distances
and discrepancies for $D_N=X_N=\{x_1,\dots, x_N\}\in \Omega_N$  as random variables
on the probability space $\Omega_N$.  We put
\begin{align}
\rho[X_N]  & = \sum\limits_{i\ne j} \rho (x_i,x_j),
\tag{3.9}
\\
\rho^*_r[X_N]  & = \sum\limits_{i\ne j} \rho^* (x_i,x_j),
\tag{3.10}
\\
\rho^*[\xi, X_N]  & = \sum\limits_{i\ne j} \rho^* (\xi, x_i,x_j)
\tag{3.11}
\end{align}
and 
\begin{align}
\la_r[X_N]& = \sum\limits_i \la_r (x_i,x_i) +
\sum\limits_{i\ne j} \la_r (x_i,x_j),
\tag{3.12}
\\
\la[\xi, X_N]& = \sum\limits_i \la (\xi, x_i,x_i) +
\sum\limits_{i\ne j} \la (\xi, x_i,x_j).
\tag{3.13}
\end{align}

The probabilistic invariance principle can be stated as follows.

\begin{proclaim}{Theorem 3.1} Let   $\R_N$ be an equal measure partition of
a compact metric space $\M$. Then the expectations of the random variables
{\rm (3.10), (3.11), (3.12)} and {\rm(3.13)} on the probability space $\Omega_N$ 
satisfy the following relations 
\begin{align}
& 2\mathbb E_N\la_r [\, . \,] + \mathbb E_N\rho^*_r[\, \cdot \,] =
N^2 \langle \rho^*_r\rangle ,
\tag{3.14}
\\
& 2\mathbb E_N\la [\xi,\, . \,] +\mathbb E_N\rho^* [\xi, \cdot]   =
N^2\langle \rho^* (\xi)\rangle.  
\tag{3.15}
\end{align}
\end{proclaim}

\begin{proof}[Proof] Summing the equality (2.16) over $x_1,x_2\in X_N$, we obtain
\begin{equation}
2\la_r[X_N]+\rho^*_r[X_N] =N^2A^{(0)}_r +2A^{(1)}_r [X_N],
\tag{3.16}
\end{equation}
where 
$$
A^{(1)}_r[X_N] =\sum\limits_{i}A^{(1)}_r(x_i).
$$
Now, we calculate the expectation $\mathbb E_N$ of both sides in the equality 
(3.16). Using the equality(3.7) and the formulas (2.17), (2.18) и (2.6), 
we find that 
\begin{align*}
& 2\mathbb E_N\la_r[\,\cdot \,] + \mathbb E_N \rho^*_r[\,\cdot\,] 
\\
& = 
N^2A^{(0)}_r+ 2\mathbb E_NA^{(1)}_r [\,\cdot \,]
= N^2A^{(0)}_r+ 2N^2 \int\limits_{\M} A^{(1)}_r (y) \, d\mu (y) 
\\
& =
2N^2\int\limits_{\M}  \mu (B_r(y)) \, d\mu (y)
+2N^2\int\limits_{\M} \mu (B_r(y))\, d\mu (y) -4N^2 \int\limits_{\M} 
(B_r(y))^2 \, d\mu (y)
\\
& = 2N^2\int\limits_{\M} [\mu (B_r(y))-\mu (B_r(y))^2] \, d\mu (y)
=  \langle \rho^*_r\rangle .
\end{align*}
This proves the relation (3.14). 

Integrating the relation (3.14) over $r\in T$ with the measure $\xi$, we obtain
the relation (3.15). 

The proof of Theorem  3.1 is complete. 
\end{proof}

Distributions $X_N\in \Omega_N$ form a subset in the set of all $N$-point 
distributions $D_N\subset \M$. Therefore, 
\begin{align}
\rho_N(\M) & \ge \mathbb E_N\rho[\, \cdot \,],
\tag{3.17}
\\
\la_N(\xi,\M) &  \le \mathbb E_N\la[\xi,\cdot]. 
\tag{3.18}
\end{align}
These inequalities in conjunction with Lemma 3.1 and Theorem 3.1 lead to
the following basic bounds.

\begin{proclaim}{Theorem 3.2} Let $\R$ be an equal measure partition of
a compact metric space $\M$. Then we have the following bound
\begin{equation}
\rho_N(M)\ge N^2\langle \rho \rangle -N \| \R_N\|_1.
\tag{3.19}
\end{equation}
If additionally, the measure $\xi$ satisfies the condition {\rm(2.10)}, then
we also have the following bound 
\begin{equation}
\la_N(\xi,\M) \le \frac 12 c_0(\xi) N \| \R_N\|_1 .
\tag{3.20}
\end{equation}
\end{proclaim}

\begin{proof} [Proof] Applying the formula (3.8) to the random variable (3.9), 
we obtain 
$$
\mathbb E_N\rho [\, \cdot \,] =N^2\langle \rho \rangle -N^2 Q_N(\rho), 
$$
where 
$$  
Q_N(\rho) \! =\! \sum\limits_i \iint\limits_{V_i\times V_i} \rho
(y_1, y_2)\, d\mu (y_1)  \, d\mu (y_2)
 \le N^{-2} \sum\limits_{i} \diam (\rho, V_i)\!=\!N^{-1} \| \R_N\|_1.
$$ 
Therefore,
$$
\mathbb E_N\rho [\,  . \,] \ge N^2 \langle \rho \rangle -N \| \R_N\|_1.
$$
Comparing this bound with the inequality (3.17), we obtain the bound (3.19). 

Let the measure $\xi$ satisfy the condition (2.10). Applying the formula (3.8) 
to the random variable (3.11), we obtain 
$$
\mathbb E_N\rho^*[\xi, \cdot\,] =N^2 \langle \rho (\xi)\rangle -N^2Q_N
(\rho^*(\xi)),
$$
where 
$$
Q_N(\rho^*(\xi)) =\sum\limits_i \iint\limits_{V_i\times V_i} \rho^*
(\xi, y_1,y_2) \, d\mu (y_1) \, d\mu (y_2)
$$
$$
\le N^{-2} \sum_i \diam (\rho^*(\xi), V_i)\le N^{-2}c_0(\xi) \sum\limits_i
\diam (\rho, V_i)
=N^{-1}c_0(\xi) \|\R\|_1. 
$$
Therefore, 
$$
\mathbb E_N \rho^* [\xi, \cdot \,]  \ge N^2 \langle \rho^*(\xi)\rangle -N
c_0(\xi)\|\R_N\|_1.
$$
Substituting this bound to the equality (3.15), we obtain 
$$
2\mathbb E_N\la [\xi,\cdot \,] \le  Nc_0(\xi) \|\R_N\|_1 
$$
Comparing this bound with the inequality (3.18), we obtain the bound (3.20).

The proof of Theorem 3.2 is complete.
\end{proof}

\section*{4.Construction of equal measure partitions. Proof of Theorem 1.1}

In this section we will prove the following general theorem.

\begin{proclaim}{Theorem 4.1} Let $\M$ be a compact $d$-rectifiable metric space.
Then, for each $N$ there exists an equal measure partition $\R_N$ of the space $\M$, such that
\begin{equation}
\| \R_N\|_1 \le d2^{d-1} \Lip (f) N^{-\frac 1d},
\tag{4.1}
\end{equation} 
where $\Lip (f)$ is the Lipschitz constant of 
the map $f$ in the definition of $d$-rectifiability of the space $\M$.
\end{proclaim}

Theorem 1.1 follows immediately from Theorems 3.2 and 4.1. 

\begin{proof}[Proof of Theorem 1.1] It suffices to substitute the boundД (4.1) 
to the bounds (3.19) and (3.20).
\end{proof}

At the present time, for spheres $S^d$ equal measure partitions are constructed
to satisfy the bound
 \begin{equation}
\| \R_N\|_{\infty}  \le c(d) N^{-\frac 1d}
\tag{4.2}
\end{equation}
with a constant $c(d)$ independent of $N$. For subsequences  $N=c_dm^d$, where $m>0$ are integers, 
such partitions were described still in the paper \cite{1} by Alexander. In the general case of all sufficiently large $N$  such partitions for spheres $S^d$ were constructed in the paper \cite{12} by
Rakhmanov, Saff and Zhou. 
 
Certainly, the bound (4.2) is stronger than (4.1), because  
$\| \R_N\|_1\le \| \R_n\|_{\infty}$, see (3.2) and (3.3). 
However, the corresponding constructions in  \cite{1, 12} significantly depend 
on the geometry of spheres $S^d$ as smooth submanifolds in $\mathbb R^{d+1}$,
while the bound (4.1) can be established for arbitrary compact $d$-rectifiable 
metric spaces. 

The proof of Theorem 4.1 is relying on three lemmas. Lemma 4.1 contains a very simple 
result which is needed at each step of our inductive construction. 
Our construction of partitions is described in Lemma 4.2 for a special case of 
a measure concentrated on the $d$-dimensional unite cube. The bound (4.1) for
such equal measure partitions of the unit cube is given in Lemma 4.3. 
Once these partitions of the unite cube are constructed, the proof of Theorem 4.1 
can be easily completed on the base of Definition 1.1.

Let $\nu_0$ be a finite non-negative measure on the unite segment $I=[0,1]$. 
Suppose that the measure $\nu_0$ has not a discrete component. Then, its
distribution function $\varphi (z)=\nu_0([0,z])$, $z\in I$, is continuous, non-decreasing,
$\varphi (0)=0$ and $\varphi(1)=\nu_0(I)$. Notice that there is a one-to-one correspondence 
between such functions and finite measures on $I$ without discrete components.

Since the graph of $\varphi$ can have horizontal parts, we define the inverse 
function $\varphi^{-1}$  by
\begin{equation}
\varphi^{-1}(t)=\sup \{z:\varphi (z) =t\}, \quad 
t\in [0,\nu_0(I)].
\tag{4.3}
\end{equation}
Let $1\le n\le k$ be integers and
\begin{equation}
n=\sum\limits^{k}_{i=1} n(i)
\tag{4.4}
\end{equation}
be an arbitrary representation of $n$ as a sum of $k$ non-negative summands $n(i)\ge 0$. 

Define points $\la (0)=0< \la (1) \le \dots \le \la (k) =1$, by
\begin{equation}
\la (j)= \varphi^{-1} \Big( n^{-1}\sum\limits^{j}_{ji=1} n(j) 
\nu_0(I)\Big), \quad 1\le j\le k,
\tag{4.5}
\end{equation}
and consider the segments $\Delta (j)=[\la (j-1),\la (j)]\subset I$,
of length $l(j)=\la (j)-\la (j-1)$, $1\le j\le k$.

We have immediately the following result.
 
\begin{proclaim}{Lemma 4.1} With the above assumptions, the segments
$\{\Delta (j), i\le j\le k\}$ form a partition of the unite segment  
$$
I=\bigcup\limits^{k}_{j=1} \Delta (j),  \quad 
\nu_0 (\Delta (j_1) \cap  \Delta  (j_2))=0, \quad j_1\ne j_2,
$$
furthermore  
\begin{equation}
\sum\limits^{k}_{j=1} l(j)=1
\tag{4.6}
\end{equation}
and
\begin{equation}
\nu_0 (\Delta (j))= \frac{n(j)}{n} \nu_0 (I).
\tag{4.7}
\end{equation}
\end{proclaim}

\smallskip
\noindent{\bf Remark.} If the measure $\nu_0\equiv 0$ identically, then for 
any $n\ge 1$ the partition given in Lemma 4.1 takes the form
\begin{equation}
\Delta (1)=[0,1], \quad \Delta (j)=[1,1]=\{ 1\}, \quad  2\le j\le k.
\tag{4.8}
\end{equation}
it is convenient to agree that for $\nu_0\equiv 0$  the partition (4.8) takes also place 
for $n=0$. 
\smallskip

Now we wish to generalize Lemma 4.1 to the $d$-dimensional unite cube $I^d=[0,1]^d$. 
Introduce some notation. 

Let $N\ge 1$ be an integer, $k=\lceil N^{1/d}\rceil$, $N\le k^d$, and 
\begin{equation}
N=\sum\limits^{k}_{i_1,\dots,i_d=1} \, N(i_1, \dots, i_d)
\tag{4.9}
\end{equation}
be an arbitrary representation of $N$ as a sum of numbers 
$N(i_1,\dots, i_d)$ equal to 0 or 1. 

Introduce the following non-negative integers 
\begin{equation}
N(i_1,\dots, i_q) =
\sum\limits^{k}_{i_{q+1,\dots,i_d=1}} N(i_1, \dots, i_d), \quad i\le q< d.
\tag{4.10}
\end{equation}
These integers satisfy the following relations 
\begin{equation}
N(i_1,\dots,i_q)= \sum\limits^{k}_{i_{q+1}=1}\, N(i_1,\dots,i_{q+1}),\quad
i\le q<d,
\tag{4.11}
\end{equation}
and 
\begin{equation}
N=\sum\limits_{i_1=1}^k \, N(i_1).
\tag{4.12}
\end{equation}

\begin{proclaim}{Lemma 4.2} Let $\nu$ be a finite non-negative measure on $I^d$ with
a continuous distribution function  
\begin{equation}
\varphi(Z)=\nu ([0,z_1]\times \dots \times [0,z_d]), \quad 
Z=  (z_1,\dots,z_d) \in I^d.
\tag{4.13}
\end{equation}
Then, for any representation of an integer $N$  as the sum $(4.9)$ there exists
a sequence of partitions  
$$
\mathcal P (q)=\{ \Pi (i_1, \dots, i_q), \quad 
1\le j\le k, \quad 1\le j\le q\}, \quad q=1,\dots,d,
$$
of the unite $I^d$ into rectangular boxes of the form  
\begin{align*}
& \Pi (i_1,\dots,i_q) = \prod\limits^{q}_{j=1} \Delta (i_1,\dots , i_j)
\times [0,1]^{d-q}
\\
& =\{ Z=(z_1,\dots,z_d) \in I^d:z_j\in \Delta (i_2,\dots,i_j), 
\\
& 1\le j\le q, \quad y_j\in [0,1] , q+1\le j\le d\}, 
\tag{4.14}
\end{align*} 
where $\Delta (i_1,\dots, i_j)\subset I$ are some segments. 

For any fixed indexes  $i_1,\dots, i_{j-1}$ the segments
$\{ \Delta(i_1,\dots,i_{j-1}, i_j)$, $i_j=1,\dots,k\}$ form a partition of $I$,
\begin{equation}
\sum\limits^{k}_{i_j=1} \, l(i_1,\dots, i_{j-1}, i_j) =1,
\tag{4.15}
\end{equation}
where $l(i_1,\dots,i_j)$ are lengths of the segments $\Delta (i_1,\dots,i_j)$. 

The measures of the rectangular boxes {\rm(4.14)} satisfy the relations
\begin{equation}
\nu(\Pi (i_1,\dots, i_q))  =
\frac{N(i_1,\dots, i_q)}{N} \nu (I^d).
\tag{4.16}
\end{equation}
\end{proclaim}                                     

\begin{proof}[Proof] We construct the partitions $\mathcal P (q)$, $q=1,\dots, d$, 
by induction on $q$. The partition $\mathcal P (1)$ is defined as follows.  

Consider the following one-dimensional distribution function 
$$
\varphi (z_1)=\nu ([0,z_1])\times [0,1]^{d-1})= \nu_0 ([0,z_1]), \quad 
z_1\in I, 
$$
where $\nu_0$ is the corresponding measure on $I$ с $\nu_0(I)=\nu(I^d)$. 
Applying Lemma 4.1 with  $n=N$ и $n(i_1)=N(i_1)$, see~(4.4), to this function,
we obtain a partition of $I$ into segments $\Delta (i_1)$ of length $l(i_1)$, 
moreover, 
$$
\sum\limits^{k}_{i_1=1} l(i_1) =1.
$$
see~(4.6). We put
$$
\Pi (i_1)=\Delta (i_1)  \times [0,1]^{d-1} \subset I^d,
$$
then,  
$$
\nu (\Pi(i_1)) =\nu_0 (\Delta (I_1)) = \frac{N(i_1)}{N} \nu (I^d). 
$$
see~(4.7). Thus, the partition $\mathcal P (1)$ is constructed.

Suppose that the partition $\mathcal P(q)$ is already constructed for some $q$,
$1\le q <d$. Then, the partition $\mathcal P(q+1)$ can be constructed as follows. 

For each rectangular box (4.14) we consider the following one-dimensional 
distribution function 
\begin{align*}
 \varphi(i_1,\dots, i_q,z_{q+1}) & =\nu 
\left( \prod\limits^{q}_{j=1} \Delta (i_1,\dots,i_j)\times
[0,z_{q+1}] \times [0,1]^{d-q-1} \right)
\\
& =\nu^{(i_1,\dots,i_q)}_{0} ([0,z_{q+1}]) , \quad z_{q+1}\in I,
\tag{4.17}
\end{align*}
where $\nu^{(i_1,\dots,i_q)}_{0}$ is the corresponding measure on $I$ and 
\begin{equation}
\nu^{(i_1,\dots,i_q)}_{0} (I)  =\nu (\Pi (i_1,\dots, i_q)) =
\frac{N(i_1,\dots, i_q)}{N}  \nu (I^d),
\tag{4.18}
\end{equation}
see~(4.16). 

Applying Lemma 4.1 with $n=N(i_1,\dots,i_q)$ и $n(i_{q+1})=N(i_1,\dots,i_{q+1})$,  
to the function (4.17), we obtain a partition of $I$ into segments $\Delta (i_1,\dots, i_1, i_{q+1})$ 
of length $l(i_1,\dots,i_q, i_{q+1})$, $1\le i_{q+1}\le k$, moreover, 
$$
\sum\limits^{k}_{i_{q+1}=1} l(i_1,\dots, i_q, i_{q+1})=1.
$$
We put
$$
\Pi (i_1,\dots, i_{q+1}) =\prod\limits^{q+1}_{j=1} \Delta (i_1,\dots, i_j)
\times [0,1]^{d-q-1}
$$
For these rectangular boxes, we have
\begin{equation}
\nu (\Pi (i_1,\dots,i_{q+1}))  =
\frac{N(i_1,\dots, i_{q+1})}{N} \nu (I^d).
\tag{4.19}
\end{equation}
Indeed, if $N(i_1,\dots, i_q\ge 1$, then in view of (4.18), we obtain 
$$
\nu(\Pi(i_1,\dots, i_{q+1})) =
\frac{N(i_1,\dots,i_{q+1})}{N(i_1,\dots,i_q)}  \nu (\Pi (i_1,\dots, i_q))
= \frac{N(i_1,\dots, i_{q+1})}{N} \nu (I^d). 
$$
If $N(i_1,\dots, i_q)=0$, then $\nu(\Pi(i_1,\dots,i_q))=0$ and the segments
$\Delta (i_1,\dots, i_q, i_{q+1})$, $1\le i_{q+1}\le k$, are defined by (4.8). 
Therefore, $\nu (\Pi(i_1,\dots, i_{q+1}))=0$ and 
the equality (4.19) is also true. 

Thus, the partition $\mathcal P (q+1)$ is constructed. 

Induction on $q$ completes the proof of Lemma 4.2.
\end{proof}

Consider the partition  $\mathcal P(d)=\{\Pi(i_1,\dots,i_d)$, $i\le i_j\le k$, $i\le j \le d\}$
constructed in Lemma 4.2. For this partition the equality (4.16) takes the form
\begin{equation}
\nu (\Pi(i_i,\dots, i_d)) =
\frac{N(i_1,\dots, i_d)}{N} \nu (I^d),
\tag{4.20}
\end{equation}
where, by definition, the numbers $N(i_1,\dots,i_d)$ are equal to 0 or 1, moreover, 
the number of $N(i_1,\dots, i_d)=1$ is $N$, see (4.9). 

We put 
$
\mathcal A=\{ \alpha =(i_1,\dots, i_d): N(i_1,\dots, i_d)=1\}
$,  $\# \mathcal A=N$, and introduce the following partition of the unite cube 
$\mathcal P_N=\{\Pi(\alpha),\alpha \in \mathcal A\}$. We write 
$\|\mathcal P_N\|_1$ for the average diameter of the rectangular boxesс $\Pi(\alpha)$,
$\alpha \in \mathcal A$, with respect to Euclidean metric 
$\| \cdot \|$ в $\mathbb R^d$, see also (4.2), 
\begin{equation}
\|\mathcal P_N\|_1 =\frac 1N \sum\limits_{\alpha \in \mathcal A}
\diam (\| \cdot \|, \Pi_{\alpha}).
\tag{4.21}
\end{equation}

\begin{proclaim}{Lemma 4.3} $\mathcal P_N =\{\Pi(\alpha), \alpha\in \mathcal A\}$ is an equal 
measure partition,
\begin{equation}
\nu (\Pi (\alpha))  =N^{-1} \nu (I^d)
\tag{4.22}
\end{equation}
and we have the bound
\begin{equation}
\|\mathcal P_N\|_1 < d2^{d-1} \, N^{-\frac{1}{d}}
\tag{4.23}
\end{equation}
\end{proclaim}

\begin{proof} [Proof] The equality (4.22) follows from (4.20)  and the definition of 
the partition $\mathcal P_N$.

The Euclidean diameter of a rectangular box  $\Pi(i_1,\dots, i_d)$ does not exceed 
the sum of lengths of its sides:
\begin{equation}
\diam (\| \cdot \|, \Pi (i_1,\dots, i_d))  \le l (i_1)+l(i_1,i_2) +\dots
+ l(i_1,\dots,i_d),
\tag{4.24}
\end{equation}
where $l(i_1,\dots, i_j)$ are lengths of the segments $\Delta (i_1,\dots, i_j)$,
$1\le j\le d$, see (4.14). 

Using (4.21), (4.24) and (4.15),  we obtain 
\begin{multline*}
N\|\mathcal P_N\|_1\le  \sum\limits^{k}_{i_1,\dots, i_d} \diam
(\|\cdot\|, \Pi (i_1,\dots, i_d))
\\
\le \sum\limits^{d}_{j=1} \sum\limits^{k}_{i_i,\dots,i_d=1}
l(i_1,\dots,i_j) =  \sum\limits^{d}_{j=1}   k^{d-j}
\sum\limits^{k}_{i_1,\dots,i_{j=1}}   l(i_1,\dots, i_j)
\\
= \sum\limits^{d}_{j=1} k^{d-1} = dk^{d-1}. 
\end{multline*}
Since $k=\lceil N^{1/d}\rceil$ and $k<N^{1/d}+1$, we have
$$
N\|\mathcal P_N\|_1 <d (N^{1/d}+1)^{d-1} = d(1+N^{-1/d})^{d-1} N^{1-\frac 1d}
\le d2^{d-1} N^{1-\frac 1d}, 
$$
that is equivalent to (4.23).

The proof of Lemma 4.3 is complete. 
\end{proof}

\begin{proof}[{\bf Proof of Theorem 4.1}] Let $\M$ $d$-rectifiable space. 
Without loss of generality, we can assume that in Definition 1.1 the measure $\nu$ 
is concentrated on the subset $\mathcal O \subset I^d$ and the measures $\mu$  и $\nu$
are normalized by $\mu (\mathcal M)=\nu (\mathcal O)=1$. 

Since the measure $\nu$ is absolutely continuous with respect to Lebesgue measure on $I^d$, 
its distribution function (4.13) is continuous, and Lemmas 4.2 and 4.3 can be applied.  

Let $\mathcal P_N=\{\Pi(\alpha),\alpha \in \mathcal A\}$ be an equal measure partition
of the unit cube $I^d$ given in Lemma 4.3 for the measure $\nu$. 
Consider the following collection of subsets in the space $\M$
\begin{equation}
\R_N =\{ V(\alpha),\alpha\in \mathcal A\}, \quad 
V(\alpha) = f(\Pi (\alpha) \cap \mathcal O). 
\tag{4.25}
\end{equation} 
Using the formula (1.13), we obtain
$$
\mu (V(\alpha)) =\nu (\Pi (\alpha) \cap \mathcal O)= \nu (\Pi (\alpha)) =N^{-1},
$$
since the measure  $\nu$  is concentrated on $\mathcal O$. 

By definition, the map $f:\mathcal O\to \M$  is an injection. Therefore 
$$
\mu(V(\alpha_1)\cap V(\alpha_2)) =\nu
(\Pi(\alpha_1)\cap \Pi (\alpha_2)\cap \mathcal O)=0, \quad \alpha_1\ne \alpha_2. 
$$
Thus, the collection of subsets $\R_N$  в (4.25) is an equal measure partition
of the $\M$. 

By definition, the map $f:\mathcal O\to \M$ is also a Lipschitz map, see~(1.12). 
Therefore 
$$
\diam (\rho, V(\alpha))\le \Lip (f) \diam (\| \cdot \|, \Pi (\alpha) \cap 
\mathcal O) \le \Lip (f) \diam (\| \cdot \|, \Pi (\alpha))
$$
and 
\begin{equation}
\| \R_N\|_1 \le \Lip (f) \|\mathcal P_N\|_1,
\tag{4.26}
\end{equation}
see (3.2) and (4.21). 

Substituting the bound (4.23) to (4.26), we obtain the bound (4.1). 

The proof of Theorem 4.1 is complete. 
\end{proof}


\begin{thebibliography}{3} 

\bibitem{1} 
J.~R.~Alexander,
\emph{On the sum of distances between $n$ points on a sphere}. --- Acta 
Math. Hungar.
{\bf 23} (3--4) (1972), 443--448.

\bibitem{2} 
J.~R.~Alexander, J.~Beck, W.~W.~L.~Chen,
\emph{Geometric discrepancy theory and uniform distributions}. --- in Handbook of Discrete
and Computational Geometry (J.~E.~Goodman and J.~O'Rourke eds.), Chapter 10, pages 185--207. 
CRC Press LLC, Boca Raton, FL, 1997.

\bibitem{3} 
A.~Barg, M.~Skriganov,
\emph{Association schemes on general measure spaces and zero-dimensional 
Abelian groups}. --- Advances in Math. {\bf 281} (2015), 142--247.

\bibitem{4}
J.~Beck,
\emph{Sums of distances between points on a sphere: An application of the 
theory of irregularities 
of distributions to distance geometry}, Mathematika {\bf 31} (1984), 33--41. 

\bibitem{5}
J.~Beck, W.~W.~L.~Chen,
\emph{Irregularities of Distribution}. --- Cambridge Tracts in Math., 
vol.~89, Cambridge Univ. Press, 1987.

\bibitem{6}
D.~Bilyk,
\emph{Discrepancy problems as particle interactions}, Preprint, 2014.

\bibitem{7}
J.~S.~Brauchart, J.~Dick,
\emph{A simple proof of Stolarsky's invariance principle}. --- Proc. Amer. 
Math. Soc. {\bf 141} (2013), 2085--2096.

\bibitem{8}
S.~Helgason, 
\emph{Differential Geometry, Lie Groups, and Symmetric Spaces}, Academic 
Press Inc., London, 1978.

\bibitem{9}
V.~I.~Levenshtein,
\emph{Universal bounds for codes and designs}, in Handbook of Coding 
Theory (V.~S.~Pless and W.~C.~Huffman eds.), Chapter~6, pages 499--648. 
Elsevier, Amsterdam, 1998.

\bibitem{10}
J.~Matou\v{s}ek,
\emph{Geometric Disrepancy. An Illustrated Guide}, Springer-Verlag, 
Berlin, 1999.

\bibitem{11}
P.~Mattila,
\emph{Geometry of Sets and Measures in Euclidean Spaces. Fractals and 
Rectifiability}, Cambridge Univ. Press, Cambridge, 1995.

\bibitem{12}
E.~A.~Rakhmanov, E.~B.~Saff, Y.~M.~Zhou,
\emph{Minimal discrete energy on the sphere}. --- Mathematical Research 
Letters {\bf 1} (1994), 647--662.

\bibitem{13}
M.~M.~Skriganov,
\emph{Point distributions in compact metric spaces}, II (in preparation).

\bibitem{14}
K.~B.~Stolarsky,
\emph{Sums of distances between points on a sphere}, II. Proc. Amer. Math. Soc.
{\bf 41} (1973), 575--582. 

\end{thebibliography}
\end{document}